\documentclass[11pt]{article}
\usepackage{xcolor}
\usepackage{amsmath}
\usepackage{amssymb}
\usepackage{graphicx}
\usepackage{soul}
\usepackage{pgf}
\usepackage{hyperref}
\usepackage{mathrsfs}

\newtheorem{theorem}{Theorem}[section]
\newtheorem{lemma}{Lemma}[theorem]

\newtheorem{cor}{Corollary}[theorem]
\newtheorem{rem}{Remark}[section]
\newtheorem{example}{Example}[section]
\newtheorem{remex}{Remark}[example]

\newcommand{\Frac}{\displaystyle\frac}            
\newcommand{\Int}{\displaystyle\int}              
\newcommand{\Limsup}{\displaystyle\limsup}        
\newcommand{\Sum}{\displaystyle\sum}              
\newcommand{\refrm}[1]{{\rm(\ref{#1})}}           
\newcommand{\vv}{\vspace{7pt}                     

}
\newcommand{\vspandex}{\rule[-18pt]{0pt}{36pt}}   
\newcommand{\SG}{{S}}                             
\newcommand{\symS}{\wt{\SG}}                      
\newcommand{\ov}[1]{\overline{#1}}
\newcommand{\wt}[1]{\widetilde{#1}}
\newcommand{\avef}[1]{{\left\langle #1 \right\rangle}}
\newcommand{\footer}[1]{{\def\thefootnote{}\footnotetext{#1}}}
\newcommand{\equ}{\,\Leftrightarrow\,}
\newcommand{\imp}{\,\Rightarrow\,}
\newcommand{\be}{\begin{equation}}
\newcommand{\ee}{\end{equation}}
\newcommand{\Proof}{\noindent {\bf Proof.~}}
\newcommand{\QED}{\hfill \rule{2mm}{2mm}         

\vspace{7pt}                                        

\noindent }
\newcommand{\set}[1]{\left\{#1\right\} }             
\newcommand{\ex}[1]{\,{\rm exp}\set{#1}\,}           
\newcommand{\df}{:=}

\newcommand{\toD}{\,\stackrel{\mbox{\footnotesize {\cal D}}}{\longrightarrow}\,}
\newcommand\cip{\,\stackrel{\mbox{\footnotesize {{\sf P}}}}{\longrightarrow}\,}
\newcommand{\tss}[1]{\textsuperscript{#1}}           

\def\ca#1{{\cal #1}}
\newcommand{\E}{{\mathbb E}}                         
\newcommand{\PR}{{\mathbb P}}                        
\newcommand{\N}{{\mathbb N}}
\newcommand{\R}{{\mathbb R}}
\newcommand{\BT}{{\mathbb T}}
\newcommand{\Var}{\text{\sf Var}}
\newcommand{\cT}{\mathscr{T}}
\newcommand{\cW}{\mathscr{W}}                         
\newcommand{\wek}[1]{\textsf{$\textbf{#1}$}}          
\newcommand{\red}[1]{{\color{red} #1}}
\newcommand{\I}[1]{1\hspace{-3pt} {\rm I}_{#1}}       
\newcommand{\Ib}[1]{1\hspace{-3pt} {\rm I}_{\{#1\}}}  

\begin{document}
\hfill \today
\vv
\centerline{\Large\bf Orderly divergence of L\'evy Gamma integrals}

\begin{center}
Jerzy Szulga{$^\dagger$}\footer
{\begin{tabular}{l}
$^\dagger$ Department of Mathematics and Statistics, Auburn University, USA\\
{~} MSC 2020: Primary 60F15, 26A45. Secondary 60G51, 60F25, 60H05.\\
{~} Key words and phrases: Gamma process, weighted  limit theorems, monotonicity
\end{tabular}}

\vv
\vv
\begin{minipage}{0.9\linewidth}
{\small {\bf Abstract.} ``Orderly divergence'' deals with  limit theorems for weighted stochastic Gamma integrals of otherwise 
nonintegrable functions. Although for monotonic functions this category usually coincides with the classical notion of weighted limit theorems for sums of i.i.d.\ random variables but there are exceptions and the lack of monotonicity reveals new aspects that are absent in the discrete case.}
\end{minipage}
\end{center}
{\footnotesize \tableofcontents}
\sf
\section{Orderly divergence}
\subsection{Motivation}
\subsubsection{Framework}
The positive Gamma process $\SG _t$  on $\R_+$   may be perceived as a continuous stream of fractional signals, augmenting discrete signals associated with a Poisson process. One may consider signals  as random points in an atomless $\sigma$-finite measure space $(\BT,\cT,\tau)$ of infinite measure. In general (cf.\ \cite{KalS,Szu:gam}), a pure jump L\'evy positive process  yields stochastic integrals $\SG f=\int_{\BT} f\,d\SG $ for positive function with the Laplace transform $\Phi(f)\df-\ln \E e^{-\SG f}= \int_{\BT} \phi(f)\,d\tau$ with $\phi$ specified below. The Borel isomorphism allows the reduction to $\BT=\R_+$ and then the integrals are 
 definable pathwise (not stochastically) over the complete modular metric space
\[
 L^\phi=\set{f: \Phi(|f|)<\infty}.
\] 
For a Gamma process, $\phi(\theta)=-\ln \E^{-\theta \SG _1} =
{\int_0^\infty \left(1-e^{-\theta x}\right) \,dx)} =\ln(1+\theta)$, so the modular $\Phi$ itself becomes a metric due to the concavity of $\phi$. 
In contrast, the symmetrized Wiener-type integral $\symS f=\SG f-\SG^* f$, where $\SG^*$ is an independent copy of $\SG$, 
may be not path-definable. 
\vv
Since the space and most of probabilistic properties are invariant with respect to the change of the intensity (scale parameter) thus  w.l.o.g.\ we may assume that $\SG _1$ is an exponential random variable with mean 1. We may modify the intensity at will if desired.
We denote  $\lambda f= \E(\SG f)=\int_0^\infty f(x)\,dx$ and $\lambda_t f= \int_0^t f$, $\SG_t f=\int_0^t f\,d\SG$.
For $p>0$,  $\E(\SG f)^p<\infty$ iff $f\Ib{f\le 1}\in L^1$ and $f\Ib {f>1}\in L^p$, cf.\ \cite{Szu:gam}.  
\vv
We  say ``$f$ is positive'' when $f\ge 0$ rather than ``$f$ is nonnegative''; similarly for decreasing or increasing functions. We consider only $t\to\infty$ (cf.\ the comment at the end of Subsection \ref{stab}) and we write $f_t\approx  g_t$ if $\lim _tf_t/g_t=1$. If the limit is a strictly positive constant, we will write $f_t\sim g_t$.
\vv
First, we compare several pertinent numerical characteristics.
We show that the stochastic integrals for monotonic integrands comply quite well with the benchmark  
of the traditional discrete weighted limit theorems. However, we reverse the classical objective of deriving a maximal class of probability distributions given specific weights. The framework of Gamma distributions prompts a search for a maximal class of admissible weights. Next, we examine the limit behavior without monotonicity. We conclude by a discussion of strength of the assumptions, alternatives, and open questions.  We hypothesize that for increasing functions the rate of increase must be sub-exponential.

\subsubsection{Setup}
For positive integrands distributions of Gamma integrals $\SG f$ coincide \cite{Thorin,JRY} with Thorin's GGC on $\R_+$, barring degenerate distributions of point masses. The constants may emerge as simple limits through WLLN or SLLN, e.g., $\SG _t/t\to 1$ a.s. This motivates us to investigate the  {\em orderly divergence}, ``OD'' in short, defined for $\lambda_t f\to \infty$ as an analog of the WLLN:
\be\label{ord t}
\lim_{t\to\infty } R_t\toD \delta_1\quad \equ\quad R_t\cip 1,\quad\mbox{where}\quad R_t=\Frac{\SG_tf }{\lambda_t f},
\ee
The notion of OD is well defined if $f$ is locally integrable.  Using the Laplace transform, OD is equivalent to
\be\label{weak}
\lim_{t} \int_\cdot^t  \ln\left(1+s\, f(x)/\lambda_t f\right)\,dx =s,\quad s\ge 0.
\ee
We may consider various modes of convergence such as a.s.\ or  in $L^p$, $p\ge 0$:
\be\label{ODmode}
R_t\stackrel{\sf mode}{\longrightarrow} 1.
\ee
 The OD for simple function $f$ with constant weights $w_n$ on intervals formed by an arithmetic progression brings 
 back the classical WDLLN (``W'': ``weighted'', ``D'': discrete''; it could be weak or strong; cf. \cite{JOP}) for i.i.d.\ random variables:
\be\label{JOP}
R_n\df\frac{S_n}{W_n}\stackrel{\sf mode}{\longrightarrow} 1,\quad \mbox{where}\quad S_n= \sum_{k=1}^nw_k X_k,\,
W_n=\sum_{k=1}^n w_k.
\ee
Thus, the discretization of the OD, as described in Subsection \ref{Disc} and  Theorem \ref{abc:disc} in particular, enables us to tap to a rich literature on WDLLN (cf., e.g.,  \cite{CohLin,Cuz,Ete,Jaj,JOP,Sun}, to name but a sample).
Observe that the normalizing occurs prior to centering, i.e. the centering function must be equal to the denominator which shifts the focus away  from Marcinkiewicz' type LLNs  \cite[p. 243]{Loev} in which  a stabilizer $k_n=o(W_n)$ is imposed after centering.
The following numerical characteristics are essential:
\be\label{vlb}
\begin{array}{lll}
(v)& v_t\df\Frac  {\lambda_t f^2}{\lambda_t^2f}\to 0 &(\text{i.e., } \Var(R_t)\to 0 \text{ or OD in $L^2$)}\, ,\\
(\ell)& \ell_t\df\Frac{\lambda_{t-1} f  }{\lambda_t f}\to 1 &(\text{an analog of the Lo\`eve constraint \cite[p.\ 252]{Loev}}),\\
(b)& b_t\df \Frac {f_t}{\lambda_tf}\to 0 &(\text{an analog of $b_n^{(d)}=w_n/W_n\to 0$, important in \refrm{JOP}}).\\
\end{array}
\ee
Relations between the three above conditions rely on additional properties.
For bounded positive functions, OD \refrm{weak} and \refrm{vlb} are  trivially satisfied.
\begin{theorem}{~}\label{abc}
Condition \refrm{vlb}.($\ell$) is necessary for OD \refrm{ord t}. Also,
\begin{enumerate}
\item
If a locally square integrable $f_t\to\infty$  then $b_t\to 0$ $\imp$ $v_t\to 0$.
\item If $f$ is monotonic then $\ell_t\to 1$  $\equ$ $b_t\to 0$.
\end{enumerate}
\end{theorem}
\Proof First, for the sake of clarity denote  $\lambda_t=\lambda_tf$ and
$S_t=S_t f$.
By hypothesis,
\[
\left(1-\ell_t\right)
\left( \Frac{S_t-S_{t-1}}{\lambda_t-\lambda_{t-1}} - \Frac{S_{t-1}}{\lambda_{t-1}}\right)=
\Frac{S_t}{\lambda_t}  -  \Frac{  S_{t-1}  }{\lambda_{t-1}}\cip 0.
\]
Suppose by contrary that $|1-\ell_t|>\epsilon>0$ infinitely often, i.e., along an unbounded countable $K\subset \R_+$.
Then $
\left(S_{t}-S_{t-1}\right)/\left(\lambda_{t}-\lambda_{t-1}\right) \cip 1 $ along $K$.
Hence, using the ``minus log'' of the Laplace transform and the  Jensen's inequality,
\[
s=\lim_{t\in K} \int_{t-1}^{t}
\ln\left(1+\Frac{ sf(x)} {\lambda_{t}f-\lambda_{t_-1}f }\right)\,dx
\le \ln(1+s),
\]
which is false for large $s$. Hence $1-\ell_t\to 0$.
\vv
Statement 1 stems from l'H\^ospital's Rule. For Statement 2 it suffices to consider an increasing $f$.  Then $0\le 1-\ell_t\le b_t$.
Conversely, since $\lambda_{t+1} f-\lambda_t f\ge f_t$ and the function $x\mapsto x/(a+x)$ is increasing, therefore
\[
1-\ell_{t+1}= \frac{  \lambda_{t+1} f-\lambda_t f  }
{   \lambda_tf +  (\lambda_{t+1} f-\lambda_t f)  } \ge \Frac{f_t}{\lambda_t f +f_t},
\]
 Hence $1-\ell_{t+1}\ge \Frac{b_{t}}{1+b_{t}}$, i.e.,
$b_{t}\le \Frac{1-\ell_{t+1}}{\ell_{t+1}}$.
\QED

\begin{rem}{~}\label{rem:rel}

\begin{enumerate}
\item  OD also implies that
${\lambda_{t-\delta} f  }/{\lambda_t f}\to 1 $  for every $\delta>0$.

{\sf\small Indeed, Jensen's inequality applied to $ddt/\delta$  yields the contradiction
$s\le \delta \ln(1+s/\delta)$.}
\item OD and \refrm{vlb} are equivalent for a monotonic $f$.
\item
OD may enforce a ``sub-exponential growth'' of $f$.

{\sf 
\small The exponential  $e(x)=e^x$ fails \refrm{vlb}. At the same time,  $R_t e\toD \SG (1/e)$ of Thorin's class because $
\int_0^t (1+s e^x/(e^t-1))\,dx\to
\int_0^\infty \ln(1+se^{-u})\,du$.
In particular, using the variance $\sigma_t^2=\Var({\SG}_t e)\approx e^{2t}$, the CLT for the centered process $\ov{\SG}= \SG-\E \SG$ yields the non-normal limit: ${\ov{\SG}_t e} /   {\sigma_t}  \,\toD \,\ov{\SG}(1/e)$; similarly for the symmetrized process:
${\wt{\SG}_t e} /   {\sigma_t}  \,\toD \,{\wt{\SG}(1/e)}$.}
\item \label{rem:rel:per}
The match of the quantities $1-\ell_t$ and $b_t$ or $v_t$ relied on increase or boundedness of $f$. Otherwise, possible ``spikes'' may void it.

{\small \sf For example, if $f$ is a periodic function with a square-integrable base then $v_t\to 0$ yet $b_t$ may be unbounded.
Indeed, let $f$ be the periodic extension of a nonzero $g\in L^2$ on a bounded interval, say,  $[0,1]$. 
Denote $c= \|g\|_2/\|g\|_1$ and $n=\lfloor t\rfloor$. Then,
\[
\lambda_t f^2 \le (n+1)\, \|g\|_2^2 \le (n+1)\, c^2 \|g\|_1^2
\le \Frac{(n+1) c^2 }{n^2}\, \lambda_t^2f.
\]}
\item In turn, if the base $g$ is only integrable but not square integrable, then ($v$) in \refrm{vlb} fails but OD reduces to the usual LLN with i.i.d.\ copies $X_n$ of $Sg$, and thus holds a.s.\ and in $L^1$ ($S_n/n$ is a uniformly integrable reverse martingale).
\end{enumerate}
\end{rem}
\subsection{Varieties of OD}
\subsubsection{OD in $L^p$}
\begin{theorem}\label{LpL0}
For $f\ge 0$ and bounded $b_t$, OD in $L^2$ implies OD in $L^p,\,p\ge 2$.
\end{theorem}
\Proof It suffices to show  that for an even $p>2$ there is a continuous function $\psi$, $\psi(0)=0$ such that
$\E|R_t f-1|^p\le \psi(\sqrt{v_t})$.
\vv

For the sake of clarity we suppress the subscript $t$. Denoting $F_0=\lambda{\ln(1+sf)} -s $, 
the Laplace transform $L(s)=\ex{-s(R-1)}=\ex{-F_0(s)}$.
Put $ F_1(s)=F_0'(s)=\lambda \big(b/(1+sf)\big)-1$ and $F_n(s)=\lambda\big(f/(1+sf)\big)^n$ for $n\ge 2$. Then
$F_n'=-nF_{n+1}$ for $n\ge 1$. The functions $P_n$, defined by the formula  $(-1)^nL^{(n)}=P_n L$ with $P_0=1$, satisfy the recurrence
\[
P_{n+1}=-P_n'+P_n F_1,\,\,n\ge 0,\quad P(0)=1.
\]
Hence, considering a partition $n=\sum_{j=1}^n jk_j$, with a vector $\wek k=(k_1,\dots, k_n)$ of whole numbers yielding $\wek F^{\wek k}=F_1^{k_1}\cdots F_n^{k_n}$, we obtain a polynomial
\[
P_n=\sum_{\wek k} c(\wek k) \wek F^{\wek k}.
\]
Denote $m_n=F_n(0)=\lambda b^n,\,n\ge 2$ and $m_1=F_1(0)=0$. In particular, $m_2=a$. Thus, it suffices to consider the vectors $\wek m=(0,m_2,...,m_n)$ because the only nonzero terms in the above sum stem from partitions $\wek k$ such that $k_1=0$ (i.e., no ``1'' in the partition of $n$). In other words,
\[
\E(R-1)^n=P_n(0)=\sum_{k_1=0}c(\wek k) \wek m^{\wek k}.
\]
Since $b_t\le C$, then for $j\ge 2$ the moments $m_j\le C^{j-2} m_2$, Hence, using $j-2\le j$,
\[
\E(Sf-1)^n\le C^n \sum_{k_1=0}c(\wek k)  m_2^{k_j} =: \psi(a),
\]
which establishes a sought-for continuous function $\psi$ for an even $n$.
\QED
\subsubsection{Discretization}\label{Disc}
Now, we relate the orderly divergence \refrm{ord t} to the classical WDLLN.
  We choose the arithmetic step 1 for the sake of convenience although the results are equivalent for any step, i.e., for $X_k$ with the Gamma distribution of any shape parameter.
\begin{theorem}  \label{abc:disc}
  Let $f\ge 0$ be a monotonic function.
Then the {\sf OD} \refrm{ODmode} and the  WDLLN \refrm{JOP} are equivalent relative to the considered mode of convergence.
\end{theorem}

We need some auxiliary results; first, one more equivalent condition.
\begin{lemma}\label{abcd}
For an increasing or bounded nonintegrable function $f$, either of conditions OD or \ref{vlb} is equivalent to the following condition
\[
d_t\df \Frac{f_{t+1}}{\lambda_tf}\to 0.
\]
\end{lemma}
\Proof
 For an increasing $f$, we have $b_t\le d_t$ and
\[
d_t=\frac{f_{t+1}}{\lambda_tf}=\frac{f_{t+1}}{\lambda_{t+1}f-\big(\lambda_{t+1}f- \lambda_t f\big)}\le \frac{f_{t+1}}{\lambda_{t+1}f-f_{t+1}}=\frac{b_{t+1}}{1-b_{t+1}}.
\]
The convergence occurs trivially when $f$ is bounded.
\QED
Next, we examine the discrete analog of the control function $b_t$.
\begin{lemma} \label{bnd}
Let $f$ be increasing  or bounded and not integrable. Then
\[
b_t\to 0 \quad \equ \quad b_n\to 0\quad \equ \quad   b_n^{(d)}\df\frac{w_n}{W_n}\to 0.
\]
\end{lemma}
\Proof
For a bounded and non-integrable function all three condition simply hold true.
Let $f$ be increasing.  The first implication is obvious. Since $d_t\le b_t$,
\[
\frac{f_t}{\lambda_tf} \le \frac{f_{n+1}}{\lambda_{n+1}f-\left(\lambda_{n+1}f-\lambda_nf\right)}=\frac{b_{n+1}}{1-d_{n+1}}\le \frac{b_{n+1}}{1-b_{n+1}}.
\]
Since $W_n\ge \lambda_nf$, so  $b_n^{(d)}\le b_n$. Conversely,
\[
b_n \le \frac{f_n}{W_{n-1}}=\frac{f_n}{W_n-f_n}=\frac{b_n^{(d)}}{1-b_n^{(d)}}.
\]
\QED

{\bf Proof of Theorem \ref{abc:disc}.}
Let $f$ be increasing.
 W.l.o.g.\ we may and do assume that $f(0)=0$.  Denote by $R_n'$ the normalized sum for the shifted sequences
 $X'_k=X _{k+1}$ and put $n=\lfloor t\rfloor$. Then
\[
\frac{R'_{n-1}}{1+d_{n-1}}\le R_t\le \Frac{R_n}{1-b_n}+d_n X_{n+1}.
\]
Recall that $X_k=\SG_k-\SG_{k-1}$. Denote $\Delta_k=\SG_k f-\SG_{k-1}f$, $\Delta =\SG_t f -\SG_nf$, and $i_k=\lambda_kf-\lambda_{k-1}f$.
For the upper bound, the denominator
\[
\lambda_tf\ge \sum_{k=1}^n i_k\ge \sum_{k=1}^n f_{k-1}=\sum_{k=1}^n f_k-f_n.
\]
Hence we bound the ratio as follows:
\[
\Frac{\SG_t f}{\lambda_t f}=\Frac{\sum_{k=1}^n \Delta_k+\Delta}{F_t}\le
\Frac{\sum_{k=1}^n f_k X_k} {\sum_{k=1}^n f_k-f_n}+\Frac{f_{n+1}X_{n+1}}{\lambda_n f}
\]
For the upper bound we divide the numerator and denominator by $\sum_{k=1}^n f_k\ge \sum_{k=1}^n i_k =\lambda_nf$ .
For the lower bound:
\[
\Frac{\SG f_t}{\lambda_tf}=\Frac{\sum_{k=1}^n \Delta_k+\Delta}{\lambda_tf}\ge
\Frac{\sum_{k=1}^n f_{k-1}X _k} {\sum_{k=1}^n f_k}
\]
The denominator $\sum_{k=1}^n f_k=\sum_{k=1}^n f_{k-1} +f_n$. We apply the division by
$\sum_{k=1}^n f_{k-1}\ge \sum_{k=2}^n i_{k-2}=\lambda_{n-1}f$ to obtain the lower bound.
\vv
Let $f$ be decreasing. Then
\[
R_n\le R_t\le \Frac{R'_{n-1}}{1-b_n}+b_n X '_n.
\]
For the upper bound:
\[
\Frac{\SG_t f}{\lambda_tf}\le \Frac{\sum_{k=1}^n f_{k-1} X _k}{\lambda_tf}+
\Frac{f_n X _{n+1}}{\lambda_nf},
\]
we use $\lambda_tf \ge \sum_{k=1}^n f_k=\sum_{k=1}^n f_{k-1}+f_n$, then divide by
$\sum_{k=1}^n f_{k-1}\le \lambda_nf$.
For the lower bound we use the inequality $\SG _t f\ge \sum_{k=1}^n f_k X _k$ and $\lambda_t f\le \sum_{k=1}^n f_k$.
\QED
Here are some benefits of the discretization.
\begin{cor}\label{fmod} 
Let $f$ be monotonic.
 Then OD  \refrm{ODmode} in $L^p$ or a.s.\
\begin{enumerate}
\item
occurs  simultaneously for any arithmetic sequence $t_n=p+qn$;
in other words, the lag $(k-1,k)$ (see the comment after \refrm{JOP}) can be replaced by any lag $(p+q(k-1),p+qk)$;
\item occurs  simultaneously for $f$ and for any measure-preserving rearrangement of $f$ between the nodes that form an arithmetic sequence.
\end{enumerate}
\end{cor}

We obtain an alternative proof of Theorem \ref{LpL0} where a stronger assumption adds OD with respect to the a.s.\ convergence.  The hypercontraction argument couldn't be applied before discretization because Gamma random variables are not uniformly hypercontractive, which would be needed for Gamma integrals.
\begin{cor}\label{main:all p}
Let $f$ be {monotonic} with $\lambda_t f\to\infty$.
Then \refrm{ODmode} holds for any  mode of convergence, $L^p$, $0\le p< \infty$ or a.s., if and only if either condition in \refrm{vlb}  is fulfilled.
\end{cor}
\Proof
We use the hypercontractivity (cf.\ \cite{KrakSzu,Szu}, going back \cite{MarZ})).
The centered standard exponential (or $\Gamma$) random variable $Y=X-1$ is $(p,2)$-hypercontractive for $p>2$, i.e.,
$\|1+Yt\|_p \le \|1+C_pYt\|_2,\,t\in\R$.  Therefore, for linear combinations of independent copies of $Y$,
\[
\|\sum_i \alpha_i Y_i\|_p\le C_{p} \,\|\sum_i \alpha_i Y_i\|_2.
\]
Due to the monotonicity, Theorem \ref{abc:disc} carries the equivalence to \refrm{ODmode}.
\QED

\subsubsection{OD relative to the a.s.\ convergence}
Extensions of the classical Kolmogorov's SLLN typically seek to weaken the independence or equal distribution of the underlying random variables, leading to martingale differences, mixing sequences, weakly independent random variables, and so on. In addition, Marcinkiewicz type SLLNs are derived with a nomalizing sequence weaker than $W_n$, applied after centering. The same process takes place for the weighted versions of SLLN. The inverse approach, extending weights given a specific sequence $(X_k)$,   is rare. Let us remind again that in our context we normalize prior to centering.
\vv
 The classical result
 \cite[Thm.3]{JOP} (assuming again that $b_n\approx w_n/W_n$ monotonically decrease) stated that the condition
\be\label{univ}
\limsup_t t\,b_t<\infty
\ee
is necessary and sufficient for the weighted  SLLN to hold for all i.i.d.\ $X_k$ with finite mean.
The authors used the cumulative function
$N(x)=\left|\set{n: W_n/w_n\le x } \right|$,
and \refrm{univ} interprets the aforementioned condition according to our notation. That is, the inverse $N^{-1}$ can be seen as the extension of the increasing rearrangement of the discrete sequence $1/b_n$ to the continuous variable.
\vv
Essentially, only weights derived from regularly increasing functions (i.e., satisfying the so called $\Delta_2$-condition) were admitted; rapidly increasing functions were excluded from this universal criterion (cf.\ the table in \ref{test funs}). Yet, our $X_k$'s are special, i.i.d.\ Gamma, so we expect a sufficient conditions weaker than \refrm{univ}.
\vv
Within monotonicity,
the sufficient condition is equivalent to the square integrability of $b_t$, eventually. Cf., e.g.,
\cite{JOP,Jaj,Ete,Sun,CohLin,Cuz} where, essentially, this condition was assumed, often yielding Marcinkiewicz type SLLN or SLLN beyond the ``i.i.d.'' paradigm (martingale differences, mixing random variables, various types of orthogonality, etc.).
\vv
In \cite[Thm.2]{JOP} the rearrangement yielded the following sufficient condition for SLLN \refrm{JOP}:
\[
Q=\E\int_0^\infty x^2 K(|Y |)<\infty,\quad\mbox{where}\quad K(x)= \int_{x}^\infty\frac{N(y)}{y^3} \quad\mbox{and}\quad Y=X-\E X.
\]
By Fubini's Theorem, cf.\ \cite{Ete},  $Q=\int_0^1 x\, \E N(|Y |/x)\,dx$.
When the function $b_t$ is monotone decreasing then $b\approx 1/N^{-1}$, thus the sufficient condition reads
\[
\int_0^\infty b_t^2K_2(1/b_t) \,dt<\infty\quad\mbox{where}\quad K_2(x)=\int _0^x u \PR(|Y |>u)\,du.
\]
Therefore, in our case of the exponential (or Gamma) distribution, the sufficient condition reduces to the square-integrability of $b_t$, leaving many functions (cf.\ the table in \ref{test funs}) outside of its applicability.
Observe that the lack of monotonicity of $b_t$, or even mere monotonicity of the discrete sequence $b_n^{(d)}$, makes the above integrals extremely hard to estimate, if at all possible. In this case the SLLN (\refrm{ord t} a.s.) would stay undecided.
We will suppress the superscript indicating the discrete context as in $b_n^{(d)}$ since the superscript-free quantities are equivalent for monotonic functions $f$.
\begin{theorem}\label{suff:ass}
Let $f$ be increasing or bounded and nonintegrable.
Then the SLLN \refrm{ord t} or the discrete SLLN \refrm{JOP} follows from either of the following conditions:
\begin{enumerate}
\item [(i)] $\sum_n b_n^2<\infty$ or $\Int_\cdot ^\infty b_t^2\,dt<\infty$ if $b_t$ be decreasing;
\item [(ii)] $\Sum_n e^{-1/b_n}<\infty$ or $\Int_\cdot ^\infty e^{-1/b_t}\,dt<\infty$ if $b_t$ be decreasing,
\end{enumerate}
where $b_n$ or $b_t$ can be replaced by $v_n$ or $v_t$.
\end{theorem}
\Proof
Due to the decrease of $b_t$, either series or integral can be used. Since both assumptions imply that $b_t\to 0$, hence by Theorem \ref{abc} the discrete WLLN \refrm{JOP} holds.
Assuming (i), the discrete SLLN follows from the Kolmogorov Theorem and Kronecker's Lemma:
\[
\begin{array}{c}
\vspandex
\sum_n b_n^2<\infty \quad\imp\quad \sum_n b_n \ov{X}_n \mbox{ converges a.s.}\\
\quad\imp\quad \Frac{1}{W_n} \sum_{k=1}^n w_k \ov{X}_k= \frac{1}{W_n} \sum_{k=1}^n W_k b_k \ov{X}_k\stackrel{\mbox{a.s.}}{\longrightarrow} 0.
\end{array}
\]
Assume (ii). We use the approach from \cite{Cuz}, based on the truncation $Y_k=\ov{X}_k\,\Ib{|\ov{X}_k|\le D}$ and $Z_k=\ov{X}_k\,\Ib{|\ov{X}_k|> D}$.
Since $Y_k$  are $(p,2)$-hypercontractive, $p>2$, then the weighted WLLN holds simultaneously in probability and in $L^2$. That is, $v_n\to 0$.  Since
\[
\frac{\sum_i w_i^2 \ov X_i^2}{W_n^2}\le D^2\,v_n,
\]
then, by Hoeffding inequality, for every $\epsilon>0$
\[
\sum_n \PR(|R_n|>\epsilon)\le 2\sum_n\ex{-\frac{\epsilon^2}{2D^2\, v_n}}.
\]
Hence the convergence of the latter series ensures the SLLN by the Borel-Cantelli Lemma.
Regarding the tail-truncated $Z_n$,
\[
\sum_n \PR(b_nZ_n>1)\le C \sum_n e^{-1/b_n}<\infty.
\]
Hence by the Kolmogorov Three Series Theorem (reduced to One Series), $\sum_n b_k Z_k$ converges a.s. and we apply again the Kronecker Lemma. This concludes the proof.
\QED
\begin{cor}\label{dec:as}
A bounded $f$ satisfies WDSLLN  \refrm{JOP}. A monotonic $f$  satisfies OD \refrm{ord t}.
\end{cor}
\Proof
If $f$ is decreasing, then we can use the integral:
\[
\int_\cdot ^\infty b_t^2\,dt \le \|f\|_\infty \int_\cdot^\infty \frac{f}{F^2}=\int_\cdot ^\infty \frac{du}{u^2}<\infty
\]
by substitution. In general, we apply the Abel formula to the series. Let $w_0=0$ w.l.o.g.. That is, after bounding $w_k\le \|w\|_\infty$, we obtain a bound by a telescopic series:
\[
\begin{array}{ll}
\Sum_{k=1}\Frac{w_k}{W_k^2}&=\Sum_{k=1}^\infty\left( \Frac{W_k}{W_k^2}  -\Frac{W_{k-1}}{W_k^2} \right)
 =
\Sum_{k=1}^\infty W_k\left( \Frac{1}{W_k^2}  -\Frac{1}{W_{k+1}^2} \right)\\
&= \Sum_{k=1}^\infty \left(1+\Frac{W_k}{W_{k+1}}\right)\left( \Frac{1}{W_k}  -\Frac{1}{W_{k+1}} \right)
\le 2
\Sum_{k=1}^\infty \left( \Frac{1}{W_k}  -\Frac{1}{W_{k+1}} \right)<\infty.
\end{array}
\]

\section{Suppressing monotonicity}
Presently, without  monotonicity we lack tools  to relate OD to the discrete OD, e.g., enabling the a.s.\ convergence. Yet, within the metric convergence such as in $L^p$ this assumption can be weakened but with a caution (cf.\ Remark \ref{rem:rel}/\ref{rem:rel:per}).
\vv
It's convenient to use the language of classes of functions induced by the constraints: OD itself, $\ca V$: $v_t\to 0$, and $\ca B$: $b_t\to 0$. (We omit the classes induced by conditions $\ell_t\to 0$  or $d_t\to 0$.) Other restrictions may apply.

\begin{example}
Let $f\ge 0$ satisfy OD. Consider increasing sequences $s_k,\,t_k$ such that $0<s_k<t_k<s_{k+1}$. Then $f$ itself  and its non-monotonic modification $f^*$ with artificial oscillations, constructed by cuts, shifts, and ``blanks'':
\be\label{f*}
f^*(x)=\sum_k f(x-s_k)\,\I{(s_k,t_k]}(x),
\ee
belong simultaneously to the same class: OD (in $L^0$ or a.s.), $\ca V$, or $\ca B$.
\QED
\end{example}
\begin{lemma}
The classes $\ca V$, $\ca B$, and OD (in $L^0$) are closed under scaling and sums.
\end{lemma}
\Proof Within a metric convergence it suffices to utilize the pattern
\[
\frac{u+v}{U+V}= \frac{U}{U+V} \,\frac{u}{U}+\frac{V}{U+V} \,\frac{v}{V}=
\alpha \,\frac{u}{U}+\beta \,\frac{v}{V},
\]
augmented by the subsequence argument and the  compactness of $[0,1]$.
\QED
Such argument fails for OD (a.s.) (unless it is equivalent to a metric convergence).
\vv
Now we can enrich the classes by oscillations that are less trivial than \refrm{f*}.
\begin{cor}\label{f*g*}
 Let $f,g$ belong to one of the three classes: OD (in $\PR$ or a.s.), $\ca V$, or $\ca B$.
The blanks $(t_k, s_{k+1}]$ in \refrm{f*} can be filled by  $g^*$  that modifies $g$ in the same manner.
\end{cor}

Restricted to $\lambda_tf_2\to \infty$ (not implied by $\lambda_t f\to \infty$, in general), $\ca B\subset \ca V$ by the l'H\^ospitals' Rule. The reverse rule is false but the quotients involved in quantities ($v$) and ($b$) in \refrm{vlb}  are very special; yet without monotonicity it still fails.
\begin{example} $\ca B$ is a proper subclass of $\ca V$.
\end{example}
Indeed, w.l.o.g.\ we may and do assume that  $f(0)=0$ and note the tautology
\be\label{taut}
f\in \ca B\quad \iff \quad f=be^B\,\,\text{and}\,\,b_t\to 0
\ee
stemming from calculus: $\lambda_t(be^B)=\lambda(e^B)'=e^B$; i.e., $F=e^B$.  Thus, we may begin with a nonintegrable function $b$ rather than with $f$. Let
\[
g=\sum _{k=1}^\infty a_k I{(s_k,t_k]},
\quad 0<s_k<t_k<s_{k+1},\,\,m_k=t_k-s_k.
\]
Also, $(am)_n\df\sum_{k=1}^\infty a_km_k\to \infty$ (we put $(am)_0=0$).
Put $b=g/2$. When $m_k\ge 1$ then the corresponding $f\not\in \ca B$.  Meanwhile, condition (a) can be rewritten as
\[
4v_t =\Frac{4\lambda_t b^2e^{2B}}{e^{2B_t}}
=\Frac{\lambda_t g^2e^{G}}{e^{G_t}}\to 0.
\]
Consequently, $G_t= (am)_{n-1}+a_n(t-s_n)$. It suffices to consider the sequence $t_n\to\infty$ in lieu of $t\to\infty$. So,  $e^{G_{t_n}}=e^{ (am)_n}$ and
\[
\lambda_{t_n} g^2 e^{G}=\Sum_{k=1}^n {a_k}  \left(e^{(am)_k}-e^{(am)_{k-1} } \right).
\]
We assume that $a_k>0$. This assumption is necessary, because we may consider an arbitrary strictly positive sequence:
\[
c_k= e^{(am)_k}-e^{(am)_{k-1} },
\]
since $m_k$ can be chosen at will. Therefore, it suffices to falsify the implication
\be\label{vc}
\frac{\sum_{k=1}^n a_k c_k}{\sum_{k=1}^n c_k} \quad\imp\quad m_n\to 0
\ee
where $c_k$ is not summable. Next, 
let $K=\{i_k\}\subset \N$ consist of a rapidly increasing sequences such that $\{v_kc_k:k\in K\}$ is summable. For example, $i_k=2^k$ while $a_k=1$ and $c_k=1$ or $c_k=\frac 1 k$. Let
\[
a_{i_k}=1\quad\text{and}\quad \sum_{j\notin K}  a_jc_j<\infty,
\]
so the latter sum's  corresponding portion is irrelevant in \refrm{vc}. Finally, choose $a_k$ on $\N\setminus K$ such that
\[
\sum_{k\notin K} a_k c_k<\infty.
\]
Hence, the denominator in \refrm{vc} converges to a strictly positive constant.
\QED
\begin{remex}{~}\label{Hh}

\begin{enumerate}
\item The construction with strictly positive $a_k$ also illustrates that blanks in \refrm{f*} cannot be filled in arbitrarily, in contrast to a patterned filling such as in Corollary \ref{f*g*}.
\item For a differentiable $f_t\to\infty$, tautologically, $(\ln f)'\to 0$ iff $f=e^H$ with $h_t\to 0$, where $H_t=\int_{\cdot}^t h$. This defines an easily verifiable subclass $\ca H\subsetneqq \ca B$, in contrast to a possibly tedious verification of the limit $b_t\to 0$.
    \QED
\end{enumerate}
\end{remex}

We call a product $fg$  with positive factors {\em regular quasi-periodic} if $g\ge 0$ is periodic with the essential supremum $\|g\|_\infty=1$, so $f$ may be called the {\em amplitude}.
\begin{theorem}
Let $fg$ be regular quasi periodic with an increasing amplitude $f\in \cW$.
Then $v_t(fg)\le C v_t(f)$, so $fg$ satisfies OD with respect to $L^p$, $p\ge 0$.
\end{theorem}
\Proof
We note that $\lambda_t(fg)^2\le \,v_t \lambda_t^2f$.  W.l.o.g.\ we may and do assume that $g$ has the base $T=[0,1]$.  We need to bound $\lambda_tf$ by $\lambda_t fg$ but in view of Remark \ref{rem:rel}/1 with $\delta=2$ it suffices to bound $\lambda_{n-1}f$, where $n=\lfloor t\rfloor$.
\vv
For $r\in (0,1)$ denote $T_0=\{t\in [0,1]: g(t)> r\}$ and put $\rho=|T_0|$. Then we partition the remainder $T\setminus T_0$ into a union of $d$ disjoint sets $T_1,\dots T_d$ with $|T_i|\le \rho$, $i=1,...,d$, where $d\ge 1$.
Then, for $i=0,...,d$ and $k=0,...,n-2$,
\[
\int_{T_i+k} f \le
\frac{1}{r} \int_ {T_0 +(k+1)} fg\le \frac{1}{r} \int_ {T +(k+1) } fg,
\]
which implies that
\[
\lambda_{n-1} f =\sum_{k=0}^ {n-2} \int_{T+k} f = \sum_{k=0}^ {n-2} \sum_{i=0}^d \int_{T_i+k} f
\le \frac{d+1}{r} \sum_{k=0}^ {n-2} \int_{T+(k+1)} fg
\le \frac{d+1}{r} \lambda_t fg.
\]
\QED

\appendix
\section{Discussion}\label{extensions}
\subsection{Strength of sufficient conditions in Theorem \ref{suff:ass}}
 \begin{enumerate}
 \item
 Observe that (i) $\imp$ (ii) since $e^{-1/b}\le k!\,b^k$. The table in \ref{test funs} show some functions that fail (i) while satisfy (ii). Unfortunately, some even fail (ii) although WLLN is secured.
  For example, for $f_t=\ex{t/\ln^\alpha t}$ (ii) holds when $\alpha>1$ since $1/v_n\approx \ln^\alpha n$.
    Yet, when $\alpha<1$, the WLLN takes place but we don't know whether SLLN does.

  \item In \cite[Thm.1.1]{Cuz} (ii) was further strengthened but simplified at the same time, replaced by the condition $\limsup_n v_n \ln n=0$.

  \item  In deriving (ii) as a sufficient condition for SLLN we utilized the exponential decrease of the tail. When the tail's drops only at a power rate, say, $\E |X|^p<\infty$, $p>1$, with no information on higher order integrability, e.g., with infinite exponential moments,  then \cite[Thm.1.1]{Cuz} in our notation gives
      \begin{enumerate}
      \item[(iii)] $ \Limsup_t\Frac{ t^q \int_\cdot ^t f^q}{F_t^q}<\infty$,
      \end{enumerate}
      where $1/q+1/p=1$, as a sufficient condition for the SLLN. The paper contains more results along with a rich bibliography. However, they shed little additional light on the specific case of the weighted SLLN for exponential random variables:
      \begin{itemize}
      \item Does the WLLN imply SLLN?
      \item What is the sufficient and neccessary condition for the SLLN?
      \end{itemize}
      The questions remain open for other non-exponential (or non-Gamma) distributions.
  \end{enumerate}

\subsection{The role of monotonicity}
A function $f\in L^\phi$ admits its decreasing rearrangement $f^\downarrow\in L^\phi$. In contrast, for $f\notin L^\phi$,  neither $f^\downarrow$ nor the increasing rearrangement $f^\uparrow$  might exist. Although they would be definable locally on each truncated domain but the rearrangements would vary as domains expand.
\vv
Yet, the monotonicity between the nodes is crucial for discretization which may fail unless local oscillations are restrained (cf. Section \ref{extensions}). Discretization allows us to tap into the ample area of classical LLNs for weighted sums of i.i.d.\ random variables. However, since the Kronecker's Lemma is the leading tool in that area thus results are usually restricted to monotonic weights.
We will examine various sufficient conditions, searching for the weakest ones. They still leave many functions satisfying WLLN without known  SLLN.
\vv
Derivation of the a.s.\ \refrm{ord t} from the discrete SLLN relies on monotonicity of the control quantities such as $b_t$. A priori,  the quantities may be not monotonic.  However, the requirement yields specific properties of the function $f_t$. For example, the function $\ell_t$ is decreasing iff, by calculus,
\[
B(s)+B(t-1)-B(s-1)-B(t)\le 0,\quad\mbox{where $B(t)=\ln F_t$ and } t<s.
\]
The latter condition occurs if the second order differences are negative which, in turn, is implied by concavity of $B_t=\int_\cdot ^t b$ which,  for a differentiable $f$, occurs iff $b_t$ is decreasing. Indeed, $b'\le 0$ iff
\[
\left( \ln \frac{F'}{F}\right)'\le 0,
\]
i.e., $F'/F$ is decreasing, which happens if and only if $\ln F$ is concave. In other words, denoting $B_t=\int_0^t b$ with a decreasing function $b_t$,
 \be\label{fb}
 F_t=e^{B_t}\quad\mbox{and}\quad f_t=b_t\,e^{B_t}.
 \ee
 By calculus, $v_t$ is decreasing iff $b_t\le 2 v_t$, which is implied by \refrm{fb}. Indeed,
 \be\label{2atbt}
2v_t= \frac{2 \int_\cdot ^t b^2e^B }{e^{B_t}}\ge  \frac{2 b_t \int_\cdot ^t b e^B }{e^{B_t}}\ge b_t.
 \ee
 \vv
 The elementary calculus, applied to the function $d_t$, yields but a tautology. However, as seen in the proof of the lemma, $b_t\approx d_t$, with $b_t\le d_t$. So, the decrease of $d_t$ makes $b_t$ ``quasi'' or ``almost'' decreasing.

\subsection{Stabilizers}\label{stab}
Although it often suffices to study just $\BT=\R_+$ with the Lebesgue measure  $|\cdot|$ and carry over obtained results  by Borel isomorphism but in general we can consider the {orderly divergence} if  there exists $k_{_T}\uparrow \infty$, called a {\em stabilizer}, such that
\be\label{ord}
\lim_{_T} \Frac{\SG  f_T}{k_{_T}}\cip 1,\quad \text{i.e.}\quad 
\lim_{T} \int_{_T} \ln\left(1+s\, f(x)/k_{_T}\right)\,dx =s,\quad s\ge 0.
\ee
A stabilizer is unique up to the equivalence $k_T\approx m_T$.
Consequently, we consider a  net of domains $T$ such that Gamma integrals $\SG_T f$ are well defined but $\SG f$ is not. In particular, $\SG  f\I T$ are unbounded in $L^0$. 
This general notation may be useful for generalizations to  either abstract or more structured domains such as $\R^d$.
\vv
Within $T=[0,\infty)$, the cumulative function $k_t=\lambda_t=\int_\cdot ^t f$ seems to be a natural stabilizer but it doesn't immediately follow from the mere occurrence of \refrm{ord}. The property alone yields only a one-sided estimate
\be\label{liminf}
1=\lim_t \E (R_t\wedge 1)\le \liminf_t \E R_t=\liminf_t {\lambda_t f/k_t}.
\ee
A denominator $k_t$ of weaker order than $\lambda_tf$ is sought-for in the Marcinkiewicz's LLN: centering first, $\ov{\SG}_t f/{k_t}=0$, so then \refrm{ord} entails $k_t\approx \lambda_t f$.
\vv
Nevertheless, the equivalence holds for bounded functions. Indeed, choose an arbitrary $\epsilon>0$.  Then $f/k_t<\epsilon$ eventually, and since $\ln(1+x)\ge \frac{\ln(1+\epsilon)}{\epsilon}\,x $ for $x\le \epsilon$, hence \refrm{ord} implies that
\[
1=\lim_t \int_\cdot^t\ln(1+f/k_t) \ge \frac{\ln(1+\epsilon)}{\epsilon} \limsup_t \frac{\avef{f}}{k_t}.
\]
Now, let $\epsilon\to 0$. Therefore, by \refrm{liminf}  $k_t\approx \lambda_tf$.
\vv
Finally, an a fortiori argument: once we derive \refrm{ord} with a stabilizer $k_t$, then it is unique up to equivalence, and we already successfully used  $k_t\df F_t$.
\vv\vv
\vv
Nets such as $(t,1]$ or $(t,\infty)$ with $t\to 0$ fall beyond our interest for the following reason, implicit in \refrm{ord}. If $\SG f(A^c)$ exists beyond  a set $A$, then in \refrm{ord} we can use the new net $T\cap A$ without changing the stabilizer. Necessarily, \refrm{weak stab} requires $|A|=\infty$.
Indeed, suppose that  $\SG  f(T\cap A) /k_{_T}\to 1$ in distribution.
Since $\ln(1+ab)\le \ln(1+a)+\ln(1+b)$, hence by \refrm{ord}
\[
s= \limsup_{_T} \int_{T\cap A} \ln\left(1+sf(x)/k_{_T}\right)\le   |A|\ln(1+s)+1
\]
for every $s>0$. By letting $s\to\infty$ we infer that
$|A|=\infty$.
\vv
 Therefore,  in the context of the positive time axis we confine to simple intervals $T=[0,t]$, replaceable by any $[a,t],\,a>0$.
\subsection{Extensions and questions}
With a minimal effort one can extends the findings in this study to a positive square integrable L\'evy process $X$ on $\R_+$, since  $\E X  f$ and $\Var(Xf)$ are proportional to $\lambda f$ and $\lambda f^2$, respectively, and the modular-yielding function $\phi(x)$ is concave. Some technical issues may arise in  scaling the variable $x\to cx$ and in dealing with higher moments in the proof of Theorem \ref{LpL0}.
\vv
Integrands with possible negative values entail obstacles due to the necessity of replacing the feasible Laplace transform by the Fourier transform, even for simple non-positive L\'evy processes such as a skewed Gamma process (cf.\ \cite{Szu:gam}).
\vv
Finally, the OD relative to the a.s.\ convergence has relied here entirely on monotonic integrands and thus the non-monotonic case remains open.
\subsection{A table of test functions}\label{test funs}
We hypothesize that for increasing functions the necessary condition for OD \refrm{ord t} or \refrm{weak} is the sub-exponential rate of increase: $f(t)=o\left(e^{-ct}\right)$ for every $c>0$.
\vv
We  list the characterizing quantities defined in \refrm{vlb} and Remark \ref{Hh}/2, for some sub-exponential test functions. The exponential itself is an upper threshold function that fails the WDLLN.  The $m$-tuple iterated logarithm is denoted by $\ln_m t$. Some near-threshold function still satisfy WDLLN but fail Cuzick's \cite{Cuz} criterion $b_t\ln t\to 0$.  \vv
There seems to be no threshold function among bounded function that satisfy WLLN, 
since non-integrability is sufficient. However, another Cuzick's \cite{Cuz} criterion $\ln t /F_t\to 0$ 
leaves the validity of SLLN an open question; e.g., $f_t=(t\ln^\alpha t)^{-1}$, $0<\alpha<1$, 
satisfies the WLLN but the SLLN is unknown to us.
\vv
Recall that, e.g., $v_t\sim b_t$ means that $v_t/b_t\to {\rm constant}>0$. The second function breaks the pattern of the quantity  $h_t$ versus $v_t$ or $b_t$.
\vv

\[
\begin{array}{l|c|c|c}
f(t)         \quad &\quad \int_\cdot ^t f^p\sim\quad &\quad  v_t\sim\,b_t\sim &h_t\sim\\ \hline
t^\alpha,\, \alpha>0 \quad& t^{p\alpha+1} &\quad t^{-1} &\quad t^{-1}   \\
e^{\ln^\alpha t},\,\alpha\le 1  & \quad t\, e^{p\ln^\alpha t}\quad
  &\quad t^{-1}&\quad\red{t^{-1}{\ln^{\alpha-1} t}}\\
  e^{\ln^\alpha t},\,\alpha\ge 1 & \quad t\, e^{p\ln^\alpha t}\,\ln ^{1-\alpha }t\quad
  &\quad t^{-1}\ln^{\alpha-1}t &\quad t^{-1} \ln ^{\alpha-1}t\\
e^{t^\alpha},\,\alpha<1    & t^{1-\alpha} e^{p\,t^\alpha}                &\quad t^{\alpha-1}
&\quad t^{\alpha-1} \\
e^{t \ln ^{-\alpha}},\, \alpha>0 &e^{p\,t \ln ^{-\alpha}t}\, \ln^\alpha t\,   &\quad \ln ^{-\alpha} t & \quad \ln ^{-\alpha} t\\
e^{t \ln_m^{-1} t},\, m\in\N &  e^{p\,t \ln_m^{-1} t}  \ln_m t &\quad \ln_m ^{-1} t & \quad \ln_m ^{-1} t\\
t^{-\alpha} e^t,\,\alpha>0& t^{-p\alpha}\, e^{pt} \quad & \quad 1\quad &\quad 1\\
e^t &\quad e^{pt} \quad & \quad 1\quad &\quad 1
\end{array}
\]
\vv

\end{document}